\documentclass[12pt]{amsart}%
\usepackage{amsfonts}
\usepackage{amsmath}
\usepackage{amssymb}
\usepackage{graphicx}%
\newtheorem{theorem}{Theorem}[section]
\newtheorem{lemma}[theorem]{Lemma}
\newtheorem{corollary}[theorem]{Corollary}
\newtheorem{proposition}[theorem]{Proposition}

\theoremstyle{definition}
\newtheorem{definition}[theorem]{Definition}

\theoremstyle{remark}

\numberwithin{equation}{section}

\begin{document}
\title{On Moments and Symmetrical Sequences}
\author{Jiten Ahuja}
\address{Jiten Ahuja\\
Department of Mathematics\\
Louisiana State University\\
Baton Rouge\\
LA 70803\\
USA}
\email{jahuja1@lsu.edu}
\author{Ricardo Estrada}
\address{R. Estrada, Department of Mathematics\\
Louisiana State University\\
Baton Rouge, LA 70803\\
U.S.A.}
\email{restrada@math.lsu.edu }
\subjclass{30D10, 30E05, 40A25, 40B05, 40C15}
\keywords{Moment problems, Symmetrical Series}

\begin{abstract}
In this article we consider questions related to the behavior of the moments
$M_{m}\left(  \left\{  z_{j}\right\}  \right)  $ when the indices are
restricted to specific subsequences of integers, such as the even or odd
moments. If $n\geq2$ we introduce the notion of symmetrical series of order
$n,$ showing that if $\left\{  z_{j}\right\}  \ $is symmetrical then
$M_{m}\left(  \left\{  z_{j}\right\}  \right)  =0$ whenever $n\nmid m;$ in
particular, the odd moments of a symmetrical series of order $2$ vanish. We
prove that when $\left\{  z_{j}\right\}  \in l^{p}$ for some $p$ then several
results characterizing the sequence from its moments hold. We show, in
particular, that if $M_{m}\left(  \left\{  z_{j}\right\}  \right)  =0$
whenever $n\nmid m$ then $\left\{  z_{j}\right\}  $ is a rearrangement of a
symmetrical series of order $n.$ We then construct examples of sequences whose moments vanish with required density. Lastly, we construct counterexamples of several of the results valid in the $l^{p}$ case
if we allow the moment series to be all \emph{conditionally convergent. }We
show that for each \emph{arbitrary} sequence of real numbers $\left\{  \mu
_{m}\right\}  _{m=0}^{\infty}$ there are real sequences $\left\{
u_{j}\right\}  _{j=0}^{\infty}$ such that%
\[
\sum_{j=0}^{\infty}u_{j}^{2m+1}=\mu_{m}\,,\ \ \ m\geq0\,.
\]

\end{abstract}
\maketitle

\section{Introduction\label{Section:Introduction}}

The behavior of the moments of a complex valued function, $\int_{X}\left(
f\left(  t\right)  \right)  ^{m}\,\mathrm{d}\lambda\left(  t\right)  ,$
$m\in\mathbb{N},$ where $\lambda$ is a measure in $X,$ has been studied
recently by several authors. Interesting results for polynomials
\cite{Markowski-Phong, MugerTuset}, rational functions \cite{Dings,
DuistermaatKallen, Mathieu}, and real analytic functions \cite{Estrada20} have
been obtained.

Of particular interest is the case of \emph{sequences,} that is, when
$X=\mathbb{N}$ and $\lambda$ is the counting measure. In this case, for a
sequence $\left\{  z_{j}\right\}  _{j=0}^{\infty},$ the moments are given as%
\begin{equation}
M_{m}=M_{m}\left(  \left\{  z_{j}\right\}  \right)  =\sum_{j=0}^{\infty}%
z_{j}^{m}\,. \label{In 1}%
\end{equation}
In 1990 Lenard \cite{Lenard} gave the construction of a sequence $\left\{
\zeta_{j}\right\}  _{j=0}^{\infty}$ all of whose moments vanish:%
\begin{equation}
\sum_{j=0}^{\infty}\zeta_{j}^{m}=0\,,\ \ \ \ m\geq1\,. \label{In 1p}%
\end{equation}
Actually Priestly \cite{Priestley}\ established that for \emph{this} sequence
we have%
\begin{equation}
\sum_{j=0}^{\infty}F\left(  \zeta_{j}\right)  =0\,, \label{In 2}%
\end{equation}
for all entire functions with $F\left(  0\right)  =0.$ Employing an infinite
dimensional version of the Levy-Steinitz rearrangement theorem \cite{Halperin,
Rosenthal} given by Katznelson and McGehee \cite{Katznelson and McGehee},
Kellinsky-Gonzalez and the second author \cite{EKev2} were able to show that
given any \emph{arbitrary} sequence of complex numbers $\left\{  \mu
_{m}\right\}  _{m=1}^{\infty}$ there are sequences $\left\{  \xi_{j}\right\}
_{j=0}^{\infty}$ such that
\begin{equation}
\sum_{j=0}^{\infty}\xi_{j}^{m}=\mu_{m}\,,\ \ \ m\geq1\,. \label{In 3}%
\end{equation}

These results show that the sequence of moments $\left\{  M_{p}\right\}
_{p=1}^{\infty}$ does not determine the sequence $\left\{  z_{j}\right\}
_{j=0}^{\infty}.$\ It is important to observe, however, that the constructions
of \cite{Lenard} and \cite{EKev2} produce \emph{conditionally convergent}
series. In fact, it was already proved in \cite{Priestley} that if $\left\{
z_{j}\right\}  ,\left\{  \xi_{j}\right\}  \in l^{p}$\ for some $p$ and if
$M_{m}\left(  \left\{  z_{j}\right\}  \right)  =M_{m}\left(  \left\{  \xi
_{j}\right\}  \right)  $ for $m\geq m_{0}$\ then each series is\ a
rearrangement of the other. In \cite{Boudabra-Markowsky} Boudabra and
Markowsky introduced a simple but powerful technique to study the behavior of
the $M_{m}$ when $\left\{  z_{j}\right\}  \in l^{p}$\ for some $p,$ showing
that%
\begin{equation}
\overline{\lim}_{m\rightarrow\infty}\left\vert M_{m}\right\vert ^{1/m}%
=\max_{j\geq0}\left\vert z_{j}\right\vert \,. \label{In 4}%
\end{equation}
An asymptotic formula for the moments, namely $M_{m}\sim c_{m}b^{m},$\ for
some constants $c_{m}$ and \ $b=\max_{j\geq0}\left\vert z_{j}\right\vert ,$
was actually obtained in \cite{EKev2}.

In this article we consider questions related to the behavior of the moments
$M_{m}\left(  \left\{  z_{j}\right\}  \right)  $ when the indices are
restricted to specific subsequences of integers, such as the even or odd
moments. If $n\geq2$ we introduce in Section \ref{Section: Symmetric series}%
\ the notion of symmetrical series of order $n,$ showing that if $\left\{
z_{j}\right\}  \ $is symmetrical then $M_{m}\left(  \left\{  z_{j}\right\}
\right)  =0$ whenever $n\nmid m;$ in particular, the odd moments of a
symmetrical series of order $2$ vanish. In Section \ref{Section: The lq case}
we prove that when $\left\{  z_{j}\right\}  \in l^{p}$ for some $p$ then
several results characterizing the sequence from its moments hold. We show, in
particular, that if $M_{m}\left(  \left\{  z_{j}\right\}  \right)  =0$
whenever $n\nmid m$ then $\left\{  z_{j}\right\}  $ is a rearrangement of a
symmetrical series of order $n.$ On the other hand, in Section
\ref{Section: Conditionally convergent moments} we show how one can construct
counterexamples of several of the results valid in the $l^{p}$ case if we
allow the moment series to be all \emph{conditionally convergent. }Among other
results we show that for each arbitrary sequence of real numbers $\left\{
\mu_{m}\right\}  _{m=0}^{\infty}$ there are real sequences $\left\{
u_{j}\right\}  _{j=0}^{\infty}$ such that
\begin{equation}
\sum_{j=0}^{\infty}u_{j}^{2m+1}=\mu_{m}\,,\ \ \ m\geq0\,. \label{In 5}%
\end{equation}
One can even find a sequence that additionally satisfies $\sum_{j=0}^{\infty
}F\left(  u_{j}\right)  =c$ where $F$ is an odd entire function and
$c\in\mathbb{R}.$

\section{Symmetrical series\label{Section: Symmetric series}}

Let $n\geq2$. We will denote the primitive $n^{\text{th}}-$root of unity,
$e^{2\pi i/n},$ as $\omega_{n}.\smallskip$

\begin{definition}
\label{Def. SS 1}A sequence of complex numbers, $\left\{  z_{j}\right\}
_{j=0}^{\infty}$, is said to be symmetrical of order $n$ if $\{z_{j}%
\}_{j=0}^{\infty}=\left\{  z_{j}\omega_{n}^{k}\right\}  _{j=0}^{\infty}$ as
sets with repetitions, for all $k$ with $0\leq k\leq n-1.\smallskip$
\end{definition}

Notice that when $n=2$ the sequence is symmetrical if after a rearrangement it
has the form $\xi_{0},$ $-\xi_{0},$ $\xi_{1},$ $-\xi_{1},$ $\xi_{2},$
$-\xi_{2},$ $\ldots.$ If the sequence $\left\{  z_{j}\right\}  _{j=0}^{\infty
}$ is actually equal to this rearrangement and $\left\vert z_{j}\right\vert
\searrow0$ then all its odd moments converge and vanish,%
\begin{equation}
\sum_{j=0}^{\infty}z_{j}^{2p+1}=0\,,\ \ p\geq0\,. \label{SS 1}%
\end{equation}
The series giving the even moments, however, could be divergent.

A similar situation occurs for a general $n.$ Since, $\omega_{n}$ satisfies
the equation $\omega^{n}-1=0,$ the sum of all complex roots of unity is zero;
in fact, the sum of $k^{\text{th}}-$powers of roots of unity is either $n$ or
$0,$ depending on whether $n\mid k$ or $n\nmid k$, respectively.$\smallskip$

\begin{lemma}
\label{Lemma SS 1}If the sequence $\left\{  z_{j}\right\}  _{j=0}^{\infty}$ is
symmetrical of order $n$ and $\left\vert z_{j}\right\vert \searrow0$ then
$M_{p}=0$ if $n\nmid p$ whenever the series converges. If the sequence has the
form $z_{nj+q}=\omega_{n}^{q}\xi_{j},$ $0\leq q<n,$ for a sequence with
$\left\vert \xi_{n}\right\vert \searrow0$ then the series giving the $M_{p}$
converge if $n\nmid p.$
\end{lemma}

\begin{proof}
Suppose that the series $\sum_{j=0}^{\infty}z_{j}^{p}$ is convergent. Then we
can show that its sum, $M_{p},$ must vanish if $n\nmid p$ as follows,%
\begin{equation}
nM_{p}=\sum_{k=0}^{n-1}M_{p}=\sum_{k=0}^{n-1}\sum_{j=0}^{\infty}z_{j}^{p}%
=\sum_{k=0}^{n-1}\sum_{j=0}^{\infty}\left(  z_{j}\omega_{n}^{k}\right)
^{p}=\sum_{j=0}^{\infty}z_{j}^{p}\left(  \sum_{k=0}^{n-1}\left(  \omega
_{n}^{p}\right)  ^{k}\right)  =0\,. \label{SS 2}%
\end{equation}
On the other hand when $z_{nj+q}=\omega_{n}^{q}\xi_{j}$ then the partial sums
$\sum_{j=0}^{nJ+Q}z_{j}^{p}$\ of the series for $J\geq0$ and $0\leq Q<n$ have
absolute values that do not exceed $\left\vert \xi_{J}\right\vert \max_{0\leq
q<n}\left\vert \sum_{k=0}^{q}\omega_{n}^{k}\right\vert ,$ and this bound goes
to zero as $J\rightarrow\infty.\smallskip$
\end{proof}

In general the series giving the moments $M_{pn}$\ of such symmetrical series
might be divergent. But we can construct examples where they converge and
satisfy $M_{pn}\neq0$ for $p\geq1.$ Take any positive real sequence $\left\{
\xi_{j}\right\}  _{j=0}^{\infty}\in l^{1}$ and define the sequence $\left\{
z_{j}\right\}  _{j=0}^{\infty}$ as $z_{nj+q}=\omega_{n}^{q}\xi_{j},$ $0\leq
q<n,$ $j\geq0.$ Then, $M_{pn}(\left\{  z_{j}\right\}  )=nM_{pn}(\left\{
\xi_{j}\right\}  )>0$ while according to the lemma $M_{r}(\left\{
z_{j}\right\}  )=0$ whenever $n\nmid r.$

\section{The $l^{p}$ case\label{Section: The lq case}}

Let us now suppose that $\left\{  z_{j}\right\}  \in l^{p}$ for some $p.$
Following the ideas of \cite{Boudabra-Markowsky} we can show that the family
of functions
\begin{equation}
f_{l,m}(\xi)=\sum_{j=0}^{\infty}\frac{(z_{j}\xi)^{l}}{1-(z_{j}\xi)^{m}}\,,
\label{SS 3}%
\end{equation}
admits the power series expansion
\begin{equation}
f_{l,m}(\xi)=\sum_{q=0}^{\infty}M_{l+qm}\xi^{l+qm}, \label{SS 4}%
\end{equation}
where $l\geq p$ and $m\geq1;$ the radius of convergence of such series is at
least $\min\limits_{j}1/\left\vert z_{j}\right\vert .$ One may prove by
various means that for $n\geq1,$%
\begin{equation}
f_{n,n}(\xi)=\frac{1}{n}\sum\limits_{j=0}^{n-1}f_{1,1}(\omega_{n}^{j}\xi)\,.
\label{SS 5}%
\end{equation}

We first provide a very simple proof of the following result of Priestley
\cite{Priestley}.$\smallskip$

\begin{proposition}
\label{Prop. lp 1}Let $\left\{  \xi_{j}\right\}  _{j=0}^{\infty}$ and
$\left\{  \eta_{j}\right\}  _{j=0}^{\infty}$ be two sequences in $l^{p},$ for
some $p\geq1.$ Suppose that $M_{m}(\left\{  \xi_{j}\right\}  )=M_{m}(\left\{
\eta_{j}\right\}  )$ eventually. Then, the non-zero terms of the two sequences
are rearrangements of each other.
\end{proposition}

\begin{proof}
Let $m_{0}\geq p$ be such that for all $m\geq m_{0},$ the moments of the two
series coincide. Then, consider the function%
\begin{equation}
f(\omega)=\sum_{j=0}^{\infty}\frac{(\xi_{j}\omega)^{m_{0}}}{1-\xi_{j}\omega
}-\sum_{j=0}^{\infty}\frac{(\eta_{j}\omega)^{m_{0}}}{1-\eta_{j}\omega}%
=\sum_{k=0}^{\infty}M_{m_{0}+k}(\left\{  \xi_{j}\right\}  )\omega^{m_{0}%
+k}-\sum_{k=0}^{\infty}M_{m_{0}+k}(\left\{  \eta_{j}\right\}  )\omega
^{m_{0}+k}. \label{lp 1}%
\end{equation}
Since $f$ vanishes in a neighborhood of $\omega=0$ it must vanish identically.
Let $\gamma\neq0$ appear $k_{\xi}>0$ times in the the sequence $\left\{
\xi_{j}\right\}  _{j=0}^{\infty}$ and $k_{\eta}\geq0$ times in $\left\{
\eta_{j}\right\}  _{j=0}^{\infty}.$ Then, $\omega=1/\gamma$ must be a
removable singularity of $f$ since the radius of convergence of the above
power series is $\infty.$ Hence, we have
\begin{equation}
\lim_{\omega\rightarrow\frac{1}{\gamma}}\left(  \omega-\frac{1}{\gamma
}\right)  f(\omega)=\lim_{\omega\rightarrow\frac{1}{\gamma}}\left(
\frac{\gamma\omega-1}{\gamma}\right)  \left[  \frac{k_{\xi}(\gamma
\omega)^{m_{0}}}{1-\gamma\omega}-\frac{k_{\eta}(\gamma\omega)^{m_{0}}%
}{1-\gamma\omega}\right]  =-\frac{1}{\gamma}\left(  k_{\xi}-k_{\eta}\right)
=0\,. \label{lp 2}%
\end{equation}
Consequently, each non-zero term in $\left\{  \xi_{j}\right\}  _{j=0}^{\infty
}$ appears the same number of times in $\left\{  \eta_{j}\right\}
_{j=0}^{\infty}$ and vice-versa. Thus, non-zero terms of the two sequences are
rearrangements of each other.$\smallskip$
\end{proof}

Actually, this result can be easily improved if we employ a theorem of Fabry
\cite[Thm. 12.6.3]{Boas} that states that if $f\left(  \omega\right)
=\sum_{n=0}^{\infty}a_{n}\omega^{n}$ is an analytic function with radius of
convergence $r$ and if the set $\left\{  n\in\mathbb{N}:a_{n}=0\right\}  $ has
density $1,$ then the circle $\left\vert z\right\vert =r$ must be a natural
boundary for $f.$ Therefore, such an $f$ cannot be a meromorphic function
unless $r=\infty$\ and it vanishes identically.We thus obtain the ensuing
result.$\smallskip$

\begin{proposition}
\label{Prop. Priestley}Let $\left\{  \xi_{j}\right\}  _{j=0}^{\infty}$ and
$\left\{  \eta_{j}\right\}  _{j=0}^{\infty}$ two sequences that belong to
$l^{p}$ for some $p.$ If%
\begin{equation}
\sum_{j=0}^{\infty}\xi_{j}^{m}=\sum_{j=0}^{\infty}\eta_{j}^{m}\,,\ \ m\in Z\,,
\label{dR 1}%
\end{equation}
and the set $Z$ has density $1,$ $\lim_{m\rightarrow\infty}\left(  1/m\right)
\left\vert \left\{  k\in Z:k<m\right\}  \right\vert =1,$ then the non-zero
terms of the two sequences are rearrangements of one another.
\end{proposition}

\begin{proof}
Indeed, Fabry's theorem yields that the meromophic function $\sum
_{k=0}^{\infty}[M_{m+k}(\left\{  \xi_{j}\right\}  )-$ $M_{m+k}(\left\{
\eta_{j}\right\}  )]\omega^{m+k}$ vanishes identically. We can then use the
same analysis as in the proof of Propostion \ref{Prop. lp 1}.$\smallskip$
\end{proof}

Interestingly, our construction of the previous section gives for any integer
$n$ sequences $\left\{  \xi_{j}\right\}  _{j=0}^{\infty}\in l^{1}$\ such that
$M_{k}=0$ whenever $n\nmid k;$ in other words, sequences whose moments vanish
in a set of density $1-1/n.$ However, when the moments vanish in \emph{some}
sets of small density, the sequence should be the zero sequence, since we
immediately obtain the following.$\smallskip$

\begin{corollary}
\label{Cor lp 1}Suppose $\left\{  \xi_{j}\right\}  _{j=0}^{\infty}\in l^{p}$
for some $p.$ If for an integer $n\geq1$ we have that $M_{qn}=0$ for all $q$
in a set of density $1,$ then, $\xi_{j}=0$ for all $j.$
\end{corollary}

\begin{proof}
We have that $\left\{  \xi_{j}^{n}\right\}  \in l^{p}.$ Then, for $q$ in a
certain set of density $1,$ $M_{q}(\left\{  \xi_{j}^{n}\right\}
)=M_{qn}(\left\{  \xi_{j}\right\}  )=0.$ Proposition \ref{Prop. Priestley}
then yields that $\xi_{j}^{n}=0$ for all $j,$ proving the required
result.$\smallskip$
\end{proof}

On the other hand, as we will presently show, in the $l^{p}$ case, for a given
$n,$ the symmetrical sequences are basically the only ones whose moments
$M_{q}$ vanish unless $n\mid q.\smallskip$

\begin{theorem}
\label{Thm. lp 1}Suppose $\left\{  \xi_{j}\right\}  _{j=0}^{\infty}\in l^{p}$
for some $p$ be any non-zero sequence such that $M_{q}=0$ whenever $n\nmid q.$
Then, the non-zero terms of $\left\{  \xi_{j}\right\}  _{j=0}^{\infty}$ must
be a rearrangement of a symmetrical sequence of order $n.$
\end{theorem}

\begin{proof}
Without loss of generality, assume that $n\mid p$. Let the complex number
$\eta\neq0$ appear $k_{0}>0$ many times in the sequence $\left\{  \xi
_{j}\right\}  _{j=0}^{\infty}.$ Similarly, suppose that for $1\leq q\leq n-1$
the term $\eta\omega_{n}^{q}$ appears $k_{q}\geq0$ times in the sequence.
Observe for each $q$ we have that
\begin{equation}
f_{p+q,n}(z)=\sum_{j=0}^{\infty}\frac{(\xi_{j}z)^{p+q}}{1-(\xi_{j}z)^{n}}%
=\sum_{j=0}^{\infty}M_{nj+p+q}z^{nj+p+k}=0\,. \label{eq:2}%
\end{equation}
Thus, if $\xi_{j}\neq0$ then $w=1/\xi_{j}$ is a removable singularity, not a
pole of $f_{p+q,n}.$ Therefore we obtain a system of equations for each $1\leq
q\leq n-1$ as follows,
\[
\lim_{z\rightarrow1/\eta}\left(  w-\frac{1}{\eta}\right)  f_{p+q,n}%
(z)=\lim_{w\rightarrow1/\eta}\sum_{j=0}^{n-1}k_{j}\left(  \frac{\eta z-1}%
{\eta}\right)  \left(  \frac{(\eta\omega_{n}^{j}z)^{p+q}}{1-(\eta\omega
_{n}^{j}z)^{n}}\right)  =0\,,
\]
or%
\begin{equation}
\lim_{z\rightarrow1/\eta}\sum_{j=0}^{n-1}\left(  \frac{-k_{j}}{\eta}\right)
\frac{(\eta\omega_{n}^{j}z)^{p+q}}{\prod\limits_{m=1}^{n-1}(1-\eta\omega
_{n}^{m}z)}=\frac{-1}{\eta\prod\limits_{m=1}^{n-1}(1-\omega_{n}^{m})}%
\sum_{j=0}^{n-1}k_{j}\omega_{n}^{jq}=0\,. \label{eq:3}%
\end{equation}
Thus, for each $1\leq q\leq n-1$ we obtain
\begin{equation}
\sum_{j=1}^{n-1}\omega_{n}^{jq}k_{j}=-k_{0}\,. \label{lp 4}%
\end{equation}
This system has a unique solution for the matrix of coefficients, $[\omega
_{n}^{jq}]_{q,j}$, is related to the Van der Monde matrix, whose determinant
is $\prod\limits_{0\leq j<q\leq n-1}(\omega_{n}^{q}-\omega_{n}^{j}),$ namely%
\[
\frac{1}{n}%
\begin{vmatrix}
1 & 1 & 1 & \cdots & 1\\
1 & \omega_{n} & \omega_{n}^{2} & \cdots & \omega_{n}^{n-1}\\
1 & \omega_{n}^{2} & \omega_{n}^{4} & \cdots & \omega_{n}^{2n-2}\\
1 & \omega_{n}^{3} & \omega_{n}^{6} & \cdots & \omega_{n}^{3n-3}\\
\vdots & \vdots & \ddots &  & \vdots\\
1 & \omega_{n}^{n-1} & \omega_{n}^{2n-2} & \cdots & \omega_{n}^{n^{2}-n}\\
&  &  &  &
\end{vmatrix}
_{n\times n}=\frac{1}{n}%
\begin{vmatrix}
n & 0 & 0 & \cdots & 0\\
1 & \omega_{n} & \omega_{n}^{2} & \cdots & \omega_{n}^{n-1}\\
1 & \omega_{n}^{2} & \omega_{n}^{4} & \cdots & \omega_{n}^{2n-2}\\
1 & \omega_{n}^{3} & \omega_{n}^{6} & \cdots & \omega_{n}^{3n-3}\\
\vdots & \vdots & \ddots &  & \vdots\\
1 & \omega_{n}^{n-1} & \omega_{n}^{2n-2} & \cdots & \omega_{n}^{n^{2}-n}\\
&  &  &  &
\end{vmatrix}
_{n\times n}%
\]%
\[
=%
\begin{vmatrix}
\omega_{n} & \omega_{n}^{2} & \omega_{n}^{3} & \cdots & \omega_{n}^{n-1}\\
\omega_{n}^{2} & \omega_{n}^{4} & \omega_{n}^{6} & \cdots & \omega_{n}%
^{2n-2}\\
\omega_{n}^{3} & \omega_{n}^{6} & \omega_{n}^{9} & \cdots & \omega_{n}%
^{3n-3}\\
\vdots &  & \ddots &  & \vdots\\
\omega_{n}^{n-1} & \omega_{n}^{2n-2} & \omega_{n}^{3n-3} & \cdots & \omega
_{n}^{n^{2}-n}\\
&  &  &  &
\end{vmatrix}
_{(n-1)\times(n-1)}=\text{det}([\omega_{n}^{jp}]_{p,j})
\]
Since $k_{0}=k_{j}$ for each $1\leq j\leq n-1$ solves (\ref{lp 4}), it is the
only solution. Hence, $\left\{  \xi_{j}\right\}  _{j=0}^{\infty}$ must be a
rearrangement of a symmetrical series of order $n.\smallskip$
\end{proof}

Notice that when $n=2$ we obtain the following contrasting results.$\smallskip
$

\begin{corollary}
\label{Cor. lp 2}Let $\{\xi_{j}\}_{j=0}^{\infty}$ be a sequence that belongs
to $l^{p}$ for some $p.$ If all its odd moments vanish then, $\{\xi_{j}%
|j\in\mathbb{N},\text{ }\xi_{j}>0\}=\{-\xi_{j}|j\in\mathbb{N},\text{ }\xi
_{j}<0\}$ as sets with repetitions.$\smallskip$
\end{corollary}

Notice that such symmetrical sequence of order $2$ could be real. Naturally if
just one even moment of a real sequence vanishes then the sequence is the zero
sequence, but even if the terms are complex we have:$\smallskip$

\begin{corollary}
\label{Cor, lp 3}If $\{\xi_{j}\}_{j=0}^{\infty}\in l^{p}$ is such that
eventually all its even moments are zero, then it must be the zero
sequence.$\smallskip$
\end{corollary}

Vastly different results hold for two sets of vanishing moments both sets with
the same density, $1/2.$

We would also like to point out that our analysis yields that for $n\geq2$ and
$\{\xi_{j}\}_{j=0}^{\infty}\in l^{p}$ for some $p,$ then the following are equivalent:

\begin{enumerate}
\item For $n\nmid q,$ $M_{q}\neq0$ infinitely often;

\item there exists some $\xi$ appearing $k_{0}>0$ times in the sequence and
some $1\leq j\leq n-1$ such that $\xi\omega_{n}^{j}$ appears $0\leq
k_{j}<k_{0}$ times in the sequence;

\item the radius of convergence of $f_{p+j,n}(z)$ is at most $1/\left\vert
\xi\right\vert $ for some $1\leq j\leq n-1.$
\end{enumerate}

\section{Sequences whose moments vanish with a given
density\label{Section: Constructing Sequences}}

We have seen that if $\left\{  \xi_{j}\right\}  \in l^{p}$\ for some $p$ then
if the moments $M_{k}\left(  \left\{  \xi_{j}\right\}  \right)  =0$ in a set
of density $1$ then $\left\{  \xi_{j}\right\}  $ is the zero sequence. On the
other hand, we have constructed sequences whose moments vanish in a set of
density $\left(  n-1\right)  /n,$ for any integer $n\geq2.$ Our aim in this
section is to construct, for any given number $0<D<1,$ sequences $\left\{
\xi_{j}\right\}  _{j=0}^{\infty}$ whose moments vanish with the density $D,$
that is%
\begin{equation}
D=\lim_{n\rightarrow\infty}\frac{\left\vert \{k\in\mathbb{N}\hspace
{0.1cm}:\hspace{0.1cm}M_{k}=0,\hspace{0.1cm}0\leq k\leq n\}\right\vert }%
{n+1}\,.
\end{equation}
We begin with a few lemmas.$\smallskip$

\begin{lemma}
\label{D 1}Let $\left\{  a_{n}\right\}  _{n=0}^{\infty}$ be a sequence
strictly decreasing to zero with $0<a_{n}<1$ for all $n$ and with $\sum
_{n=0}^{\infty}a_{n}$ divergent. Then for all $x\in\left(  0,1\right)  $ there
exist infinite subsequences $\left\{  a_{n_{k}}\right\}  _{k}$ such that%
\begin{equation}
\prod_{k=0}^{\infty}\left(  1-a_{n_{k}}\right)  =x\ . \label{L 1}%
\end{equation}

\end{lemma}

\begin{proof}
Let $n_{0}$ be the first index for which $1-a_{n_{0}}>x.$ Recursively define
$n_{q}$ as the first index for which%
\begin{equation}
\prod_{k=0}^{q}\left(  1-a_{n_{k}}\right)  >x\ .\label{L 2}%
\end{equation}
The infinite product $\prod_{k=0}^{\infty}\left(  1-a_{n_{k}}\right)  $
converges, because the partial products (\ref{L 2}) form a decreasing sequence
bounded below by $x>0.$ That $\prod_{k=0}^{\infty}\left(  1-a_{n_{k}}\right)
\geq x$ is clear. That $\prod_{k=0}^{\infty}\left(  1-a_{n_{k}}\right)  >x$ is
not possible can be seen as follows: if the product is larger than $x$ we can
find $a_{n^{\ast}}$ where $n^{\ast}$ is not any of the $n_{k}$'s such that
\[
\left(  1-a_{n^{\ast}}\right)  \prod_{k=0}^{\infty}\left(  1-a_{n_{k}}\right)
>x\ ,
\]
(because the $\sum_{k=0}^{\infty}a_{n_{k}}$ converges but $\sum_{n=m}^{\infty
}a_{n}$ diverges). If $n_{q-1}<n^{\ast}<n_{q},$ then the definition of $n_{q}$
would not be satisfied, since we should have taken $n^{\ast}$ instead of
$n_{q}$; this a contradiction.

We can do the same construction asking the $n_{q}$ not to belong to given
finite subset $F$ of $\mathbb{N},$ and there are infinite ways to choose $F,$
so infinite ways to construct the \ subsequence.$\smallskip$
\end{proof}

We thus obtain the ensuing density result.$\smallskip$

\begin{corollary}
Consider the collection, $\mathfrak{S}$, of all numbers $r$ that have the
following form:%
\begin{equation}
1-r=\prod_{j=0}^{N-1}\left(  1-\frac{1}{n_{j}}\right)  \ ,\hspace
{0.5cm}\text{where }N\in\mathbb{N}\text{, and }gcd(n_{j},n_{k})=1\text{ for
}j\neq k\ .
\end{equation}
Then, $\mathfrak{S}$ is dense in $[0,1].$
\end{corollary}

\begin{proof}
Define $\mathfrak{P}$ as the collection of all numbers, $r$, of the form
\begin{equation}
1-r=\prod_{j=0}^{N-1}\left(  1-\frac{1}{p_{j}}\right)  ,\hspace{0.5cm}%
\text{where }N\in\mathbb{N},\text{ and }p_{j}\text{ are different primes
}\text{for }j\geq0.
\end{equation}
The lemma \ref{D 1} yields that the set $\mathfrak{P}\subset\mathfrak{S}$ is
dense in $[0,1]$ and so is the set $\mathfrak{S}.\smallskip$
\end{proof}

Let us now consider any sequence $\left\{  \eta_{j}\right\}  _{j=0}^{\infty
}\in l^{1}$ whose moments never vanish and a given number, $D$ with $0<D<1.$
Then, by Lemma \ref{D 1}, there exist a strictly increasing sequence of
primes, $\{p_{k}\}_{k=0}^{\infty}$ such that $\prod_{k=0}^{\infty}\left(
1-1/p_{k}\right)  =1-D.$ Using these primes, we then define countable sets,
$\left\{  I_{k}\right\}  _{k=0}^{\infty}$ with repetitions as follows:%
\begin{gather}
I^{\prime}:=\{\eta_{j}:j\geq0\}\ ,\\
I_{k}:=\left\{  \eta_{j}\frac{\omega_{p_{k}}^{l}}{2^{k}}\hspace{0.15cm}%
:\hspace{0.15cm}\eta_{j}\in I^{\prime},\hspace{0.15cm}0\leq l\leq
p_{k}-1\right\}  \ ,\quad k\geq0\ ,\\
I=\bigcup_{k=0}^{\infty}I_{k}\ .
\end{gather}
Thus, for $k\geq0$, the sets $I_{k}$ have as all their elements, all the terms
from symmetrical sequences of order $p_{k}$. Thus, moments for each sequence
from $I_{k}$ do not vanish with denisty $1/p_{k}.$ Define $\left\{  \xi
_{j}\right\}  _{j=0}^{\infty}$ as a bijection, $\xi:\{0,1,2,3,...\}\rightarrow
I.$


\begin{theorem}
The sequence $\left\{ \xi_{j}\right\} _{j=0}^{\infty}$ thus constructed has
moments that vanish with density $1-D$.
\end{theorem}

\begin{proof}
We first note that $\left\{  \xi_{j}\right\}  _{j=0}^{\infty}\in l^{1}.$ Also,
note that $M_{q}\left(  \left\{  \xi_{j}\right\}  _{j=0}^{\infty}\right)  =0$
if and only if $p_{k}\nmid q$ for any $k\geq0.$ Then, by the
inclusion-exclusion principle, the number of moments, $M_{q},$ $q\leq n$, that
vanish is
\begin{align}
& \left(  n-\left[  \frac{n}{p_{1}}\right]  -\left[  \frac{n}{p_{2}}\right]
-...+\left[  \frac{n}{p_{1}p_{2}}\right]  +\left[  \frac{n}{p_{1}p_{3}%
}\right]  ...-\left[  \frac{n}{p_{1}p_{2}p_{3}}\right]  -\left[  \frac
{n}{p_{1}p_{2}p_{4}}\right]  ...\right)  \nonumber\\
& \sim n\prod_{k=0}^{\infty}\left(  1-\frac{1}{p_{k}}\right)  \sim n\left(1-D\right)\ ,
\end{align}
as required.
\end{proof}

\section{Conditionally convergent
moments\label{Section: Conditionally convergent moments}}

We will now show that for almost all the results of the previous section it is
possible to find counterexamples if we allow sequences $\{\xi_{j}%
\}_{j=0}^{\infty}$ that do not belong to any $l^{p}$ but with moment series
that converge \emph{conditionally.}

Our main tool is the infinite dimensional version of the Levy-Steinitz
rearrangement\ theorem \cite{Rosenthal} given by Katznelson and McGehee
\cite{Katznelson and McGehee}. The result says that if $\left\{
\mathbf{x}_{j}\right\}  _{j=0}^{\infty}$ is a sequence in the Fr\'{e}chet
space $\mathbb{R}^{\mathbb{N}}$ such that the series%
\begin{equation}
\sum_{j=0}^{\infty}\mathbf{x}_{j}\,, \label{r3}%
\end{equation}
converges, then the set $S=S\left(  \left\{  \mathbf{x}_{j}\right\}  \right)
$ of all the possible sums of convergent rearrangements of $\left\{
\mathbf{x}_{j}\right\}  _{j=0}^{\infty}$ is an affine manifold of
$\mathbb{R}^{\mathbb{N}},$ namely,%
\begin{equation}
S=\mathbf{z}+N\,, \label{r4}%
\end{equation}
where $N$ is a vector subspace of $\mathbb{R}^{\mathbb{N}}$ and $\mathbf{z}$
is any element of $S.$ In fact \cite{Halperin}, $N$ is the polar set
\begin{equation}
N=\left\{  \mathbf{y}\in\mathbb{R}^{\left(  \mathbb{N}\right)  }:\sum
_{j=0}^{\infty}\left\langle \mathbf{y,x}_{j}\right\rangle \text{ converges
absolutely}\right\}  ^{0}\,. \label{r5}%
\end{equation}
As it is the standard practice, we identify $\mathbb{R}^{\left(
\mathbb{N}\right)  },$ the space of sequences with only a finite number of
non-zero terms, with the dual space $\left(  \mathbb{R}^{\mathbb{N}}\right)
^{\prime}$ \cite{Treves}. Therefore, when $\sum_{j=0}^{\infty}\left\langle
\mathbf{y,x}_{j}\right\rangle $ does not converge absolutely for any non zero
$\mathbf{y}\in\mathbb{R}^{\left(  \mathbb{N}\right)  }$ then the set of sums
of convergent rearrangements of $\left\{  \mathbf{x}_{j}\right\}
_{j=0}^{\infty}$ is all of $\mathbb{R}^{\mathbb{N}}.$

Using these ideas, the following result was established in \cite{EKev2}%
.\smallskip

\begin{theorem}
\label{Theorem r1}Let $\left\{  \xi_{q}\right\}  _{q=0}^{\infty}$ be a
non-zero sequence of complex numbers such that the moment series for
$M_{m}\left(  \left\{  \xi_{j}\right\}  \right)  ,$%
\begin{equation}
\sum_{j=0}^{\infty}\xi_{j}^{m}\,,\ \ m=1,2,3,\ldots\,, \label{r7}%
\end{equation}
all converge but never absolutely. Then for each arbitrary sequence of complex
numbers $\left\{  \mu_{m}\right\}  _{m=1}^{\infty}$ there is a rearrangement
$\left\{  \rho_{j}\right\}  _{j=0}^{\infty}$ of the series $\left\{  \xi
_{j}\right\}  _{j=0}^{\infty}$ such that%
\begin{equation}
\sum_{j=0}^{\infty}\rho_{j}^{m}=\mu_{m}\,,\ \ \ m=1,2,3,\ldots\,. \label{r8}%
\end{equation}

\end{theorem}

The proof employs the sequence of \emph{real} sequences $\left\{
\mathbf{w}_{j}\right\}  _{j=0}^{\infty}$ of $\mathbb{R}^{\mathbb{N}}$ given
by
\begin{equation}
\mathbf{w}_{j}=\left(  \Re e\left(  \xi_{j}\right)  ,\Im m\left(  \xi
_{j}\right)  ,\Re e\left(  \xi_{j}^{2}\right)  ,\Im m\left(  \xi_{j}%
^{2}\right)  ,\Re e\left(  \xi_{j}^{3}\right)  ,\Im m\left(  \xi_{j}%
^{3}\right)  ,\ldots\right)  \,, \label{r9}%
\end{equation}
proving that for $\mathbf{y}\in\mathbb{R}^{\left(  \mathbb{N}\right)  }$ the
series $\sum_{j=0}^{\infty}\left\langle \mathbf{y,w}_{j}\right\rangle $
converges absolutely only if $\mathbf{y=0}.$

Let us first show how this result allow us to give an alternative construction
of a series whose moments all vanish.\smallskip

\begin{lemma}
\label{Lemma cc 1}Let $\omega=e^{2\pi\alpha i}$ where $\alpha$ is irrational.
Let%
\begin{equation}
\xi_{j}=\frac{\omega^{j}}{\ln\left(  j+2\right)  }\,,\ \ \ j\geq0\,.
\label{cc 1}%
\end{equation}
Then the moment series $\sum_{j=0}^{\infty}\xi_{j}^{m},$ are conditionally
convergent for all $m\geq1.$
\end{lemma}

\begin{proof}
Indeed, this is an elementary fact.\smallskip
\end{proof}

Using the sequence of this lemma and the Theorem \ref{Theorem r1} we obtain
the following.\smallskip

\begin{corollary}
\label{Cor r1}For each arbitrary sequence of complex numbers $\left\{  \mu
_{m}\right\}  _{m=1}^{\infty}$ there is a rearrangement $\left\{  \rho
_{j}\right\}  _{j=0}^{\infty}$ of the sequence $\left\{  \xi_{j}\right\}
_{j=0}^{\infty}$ given by (\ref{cc 1}) such that
\begin{equation}
\sum_{j=0}^{\infty}\rho_{j}^{m}=\mu_{m}\,, \label{r 11}%
\end{equation}
for $m=1,2,3,\ldots\,.$ In particular, there are rearrangements $\left\{
\rho_{j}\right\}  _{j=0}^{\infty}$ of the series $\left\{  \xi_{j}\right\}
_{j=0}^{\infty}$ all of whose moments vanish.
\end{corollary}

We can also construct in this fashion counterexamples of the Theorem
\ref{Thm. lp 1} when condionally convergent series are allowed.\smallskip

\begin{corollary}
\label{Cor cc 1}Let $n\geq2.$ There are sequences $\left\{  \eta_{j}\right\}
_{j=0}^{\infty}$\ such that
\begin{equation}
M_{m}\left(  \left\{  \eta_{j}\right\}  \right)  =\sum_{j=0}^{\infty}\eta
_{j}^{m}=0\,,\ \ \ \ n\nmid m\,, \label{cc 2}%
\end{equation}
whose non-zero terms are not a rearrangement of a symmetrical sequence of
order $n.$
\end{corollary}

\begin{proof}
We can take $\left\{  \eta_{j}\right\}  _{j=0}^{\infty}$\ as a rearrangement
of the sequence $\left\{  \xi_{j}\right\}  _{j=0}^{\infty}$ of the Lemma
\ref{Lemma cc 1}. Clearly $\left\{  \eta_{j}\right\}  _{j=0}^{\infty}$\ is not
a rearrangement of a symmetrical sequence of order $n.$\smallskip
\end{proof}

A similar analysis can be used to study the odd moments of \emph{real
}sequences.\smallskip

\begin{lemma}
\label{Lemma cc 2}The series giving the odd moments $\sum_{j=0}^{\infty}%
x_{j}^{2m+1},$ $m\geq0,$\ of the real sequence
\begin{equation}
x_{j}=\frac{\left(  -1\right)  ^{j}}{\ln\left(  j+2\right)  }\,,\ \ \ j\geq
0\,, \label{cc 3}%
\end{equation}
are all conditionally convergent. If $\mathbf{y}\in\mathbb{R}^{\left(
\mathbb{N}\right)  }$ does not vanish, $\mathbf{y}=\left\{  y_{m}\right\}
_{m=0}^{\infty}$ with $y_{m}=0$ for $m>M,$ then the series%
\begin{equation}
\sum_{j=0}^{\infty}\sum_{m=0}^{M}y_{m}x_{j}^{2m+1}, \label{cc 4}%
\end{equation}
is conditionally convergent.\smallskip
\end{lemma}

We thus obtain from the Levy-Steinitz rearrangement\ theorem of Katznelson and
McGehee the next result.\smallskip

\begin{proposition}
\label{Prop. cc 1}For each arbitrary sequence of real numbers $\left\{
\mu_{m}\right\}  _{m=0}^{\infty}$ there is a rearrangement $\left\{
u_{j}\right\}  _{j=0}^{\infty}$ of the sequence $\left\{  x_{j}\right\}
_{j=0}^{\infty}$ given by (\ref{cc 3}) such that
\begin{equation}
\sum_{j=0}^{\infty}u_{j}^{2m+1}=\mu_{m}\,,\ \ \ m\geq0\,. \label{cc 5}%
\end{equation}
In particular, there are real sequences $\left\{  u_{j}\right\}
_{j=0}^{\infty}$ all of whose odd moments vanish but which are not a
rearrangement of a symmetrical sequence of order $2.$\smallskip
\end{proposition}

In \cite{Priestley} Priestley established that for the sequence $\left\{
\zeta_{j}\right\}  _{j=0}^{\infty}$ constructed by Lenard \cite{Lenard}\ we
have that $\sum_{j=0}^{\infty}F\left(  \zeta_{j}\right)  =0$ for all entire
functions with $F\left(  0\right)  =0.$ As we now show, our construction can
produce very different results.\smallskip

\begin{proposition}
\label{Prop. cc 2}Let $F$ be a real odd entire function. For each sequence of
real numbers $\left\{  \mu_{m}\right\}  _{m=0}^{\infty}$ and for each
$c\in\mathbb{R}$ there exists a rearrangement $\left\{  u_{j}\right\}
_{j=0}^{\infty}$ of the sequence $\left\{  x_{j}\right\}  _{j=0}^{\infty}$
given by (\ref{cc 3}) such that (\ref{cc 5}) holds and such that $\sum
_{j=0}^{\infty}F\left(  u_{j}\right)  =c.$
\end{proposition}

\begin{proof}
Let $\mathbf{w}_{j}=\left(  F\left(  x_{j}\right)  ,x_{j},x_{j}^{3},x_{j}%
^{5},\ldots\right)  .$ We need to show that if $\mathbf{y}\in\mathbb{R}%
^{\left(  \mathbb{N}\right)  }$ the series $\sum_{j=0}^{\infty}\left\langle
\mathbf{y,w}_{j}\right\rangle $ converges absolutely only if $\mathbf{y=0}.$

Consider the series $\sum_{j=0}^{\infty}F\left(  x_{j}\right)  .$ Since
$F\left(  x\right)  \sim ax^{2k+1},$ $F^{\prime}\left(  x\right)  \sim\left(
2k+1\right)  ax^{2k}$ for some $k$ and some $a\neq0$ it follows that $F$ is
strictly increasing in a neighborhood of the origin, and consequently, for $j$
large enough the sequence $F\left(  x_{j}\right)  $ is an alternating series
with decreasing absolute values. Therefore $\sum_{j=0}^{\infty}F\left(
x_{j}\right)  $ converges. The convergence is not absolute, however, because
$\left\vert F\left(  x_{j}\right)  \right\vert \sim a\left(  \ln\left(
2j+2\right)  \right)  ^{-2k-1}$ as $j\rightarrow\infty.$ Replacing $F$ by
$G\left(  x\right)  =y_{0}F\left(  x\right)  +\sum_{q=1}^{Q}y_{q}x^{2q-1},$
which is also entire and odd, we see that the series $\sum_{j=0}^{\infty
}G\left(  x_{j}\right)  $ is conditionally convergent, but it is not
absolutely convergent unless $y_{q}=0$ for $0\leq q\leq Q;$ in other words,
$\sum_{j=0}^{\infty}\left\langle \mathbf{y,w}_{j}\right\rangle $ does nor
converges absolutely unless $\mathbf{y=0}.$\smallskip
\end{proof}

This proposition yields that we can find real sequences $\left\{
u_{j}\right\}  _{j=0}^{\infty}$ such that $\sum_{j=0}^{\infty}u_{j}^{2m+1}=0,$
for $m\geq0$ but such that $\sum_{j=0}^{\infty}F\left(  u_{j}\right)
=c\neq0.$

\end{document}